\documentclass[12pt]{article}
\usepackage[centertags]{amsmath}
\usepackage{amssymb}

\newtheorem{thm}{Theorem}[section]

\newtheorem{lem}[thm]{Lemma}

\def\qed{\hfill $\Box$}

\def\P{{\bf P}}

\def\Z{\mathbb{Z}}

\def\({\left(} \def\){\right)}

\begin{document}

\title{The Reversible Nearest Particle System \\ on a Finite Interval\footnote{Supported in part by Grant
G1999075106 from the Ministry of Science and Technology of China}}

\date{}

\author{Dayue Chen, \  Juxin Liu \ and \ Fuxi Zhang\footnote{corresponding author, E-mail: zhangfxi@math.pku.edu.cn.} \\
}

\maketitle

{\bf Abstract:} In this paper we study a one-parameter family of
attractive reversible nearest particle system on a finite
interval. As the length of the interval increases, the time that
the nearest particle system first hits the empty set increases in
different order, from logarithmic to exponential, according to the
intensity of interaction. In particular, at the critical case, the
first hitting time increases in a polynomial order.



\section{Introduction}

A nearest particle system 
on $S= \{1,2,\cdots,N\}$ is a continuous time Markov chain with
the state space $\{ A : A \subset S \}$. The jump rates are
specified as follows:
 $$
\begin{array}{ll}
q (A, A \setminus \{x\}) = 1 & {\rm if} \ x \in A; \\
q(A, A \cup \{x\}) = \beta(l_x(A), r_x(A)) & {\rm if} \ x \in S
\setminus A; \\
q (A, B) = 0 & {\rm otherwise}.
 \end{array}
 $$
Here $l_x(A)$ and $r_x(A)$ are the distances from $x$ to the
nearest points in $A$ to the left and right respectively, with
convention that $l_x(A)$ (or $r_x(A)$) is $\infty$ if $y> x$ (or
$y< x$,
respectively) for all $y\in A$. We assume that\\
1. $\beta (l, r) = \beta (r, l)$;\\
2. $\beta (l,r)$ is decreasing in $l$ and in $r$;\\
3. $\beta(\infty, \infty) = 0, \beta(l,\infty) > 0$;\\
4. $\sum_l \beta(l, \infty) < \infty$.

There are many choices of $\beta(\cdot, \cdot)$ satisfying the
above assumptions.  \\
{\bf Example 1}. (The 1-dim contact process) $\beta (1,1) = 2
\lambda$, $\beta(1,r) = \beta(l,1) = \lambda$ for $l,r >1$, and
$\beta (l,r)
= 0$ otherwise. \\
{\bf Example 2}. (The uniform birth rate) $\beta(l, r) = \lambda/(l+r-1)$.\\
{\bf Example 3}. (The reversible case)
 $$
\beta(l,r) = \lambda \frac {\psi(l)\psi(r)}{\psi(l+r)}, \ \ \
\beta(l,\infty) = \beta(\infty,l) = \lambda \psi(l),
 $$
where
\begin{equation}\label{rev}
\psi(\cdot)> 0,\ \  \sum_{n=1}^\infty \psi(n) =1, \mbox{ and }
\frac { \psi(n)}{\psi(n+1)} \searrow 1. \end{equation}
Assume further that
\begin{equation}\label{momn}
\sum_n n^2 \psi(n) < \infty.
 \end{equation}
 For example,  $\psi(n) = c n^{-\alpha}$ for
some $\alpha
> 3$ satisfies the above requirement.

It is helpful to associate a subset $A$ of $S$ with an element
$\xi$ of $\{0,1\}^S$ and use them interchangeably: $\xi(x) = 1$ if
and only if $x\in A$. Configuration $\xi$ will be given an
occupancy interpretation. We say there is a particle in $x$ if
$\xi(x)= 1$, and we say the site is vacant if $\xi(x) =0$. Then
the above transition mechanisms can be interpreted as follows:
Each particle disappears at rate 1 independently, and a particle
is born at vacant site $x$ at rate $\beta(l_x(A), r_x(A))$.

The transition mechanisms also make sense if we replace $\{1, 2,
\cdots, N\}$ with the integer lattice $\Z$. However, because of
Assumption 4, the
state space $\{0, 1\}^\Z$ consists of four disjoint parts: \\
1) all finite subsets of $\Z$;\\
2)  all subsets of $\Z$ with infinite many particles both to the
left and to the right of the origin; \\
3) all infinite subsets of $\Z$ with finite many particles to the
right of the origin;  and \\
4) all infinite subsets of $\Z$ with finite many particles to the
left of the origin.\\
The first two cases are extensively studied, and are
 called  {\it finite} and {\it infinite} nearest particle systems respectively. A comprehensive
account can be found in Chapter 7 of Liggett (1985).
 The last two cases share many properties of the first two cases,
 and are indispensable in some occasions, e.g., Lemma 4.1.

For interacting particle systems people are most concerned with
the existence of phase transition and the critical value. For the
infinite nearest particle system with the uniform birth rate
(Example 2), the critical value is 1, see Mountford (1992).  For
the reversible nearest particle system (Example 3), the critical
value is also 1. For the contact process (Example 1), the critical
value is unknown but is between 1.5 and 2, and is denoted as
$\lambda_c$ throughout this paper.

Can the critical value of an infinite model be detected by the
counterpart on a finite interval? This interplay was first
explored for the contact process in a series paper by Durrett {\it
et al}. The main results are summerized as follows. Let $\{
\zeta^N_t : t \ge 0 \}$ be the contact process on $\{ 1, 2,
\cdots, N \}$ with the parameter $\lambda$ starting from all sites
occupied, and $\tau_N$ be the first time it hits the empty set.

 \begin{thm}\label{c2}
{\rm (i)}  If $\lambda < \lambda_c$, then there is a constant
$\gamma_1 (\lambda) > 0 $ so that as $N \rightarrow \infty$, $
\tau_N / \log N \rightarrow 1 / \gamma_1 (\lambda)$ in probability
{\rm (Durrett and Liu (1988), Theorem 1)}.

{\rm (ii)} If $\lambda > \lambda_c$, then there is a constant
$\gamma_2 (\lambda) > 0 $ so that as $N \rightarrow \infty$,
$(\log \tau_N) / N  \rightarrow \gamma_2 (\lambda)$ in probability
{\rm (Durrett and Schonmann (1988), Theorem 2)}.

{\rm (iii)} If $\lambda = \lambda_c$ and $a,b \in (0, \infty)$,
then $P \( a N \le \tau_N \le b N^4 \) \rightarrow 1$ as $N
\rightarrow \infty$ {\rm (Durrett {\it et al} (1989), Theorem
1.6)}.
 \end{thm}

We believe that these statements hold for a large class of
interacting particle systems. In this paper we like to study the
asymptotical behavior of the hitting time $\sigma_N$ of the
reversible nearest particle systems (Example 3) on a finite
interval, as the length of interval increases.  The results read
as follows. Let  $\{ C_N : N \ge 1 \}$ be any sequence of
increasing numbers such that $\lim_{N \rightarrow \infty} C_N =
\infty$.

 \begin{thm} \label{shang2}  Suppose the initial state is $\{1,2,
 \cdots, N\}$.\\
(1) If $$\lambda <  \min\{1, \ \min_n  \frac
{\psi(n)}{\sum_{l+r=n}\psi(l) \psi(r)} \},$$ then $E\sigma_N \leq
C \log N$ for some constant $C$ which is independent of $N$, and
 $$
\lim_{N \rightarrow \infty} P ( \sigma_N \leq C_N \log N ) = 1;
 $$
(2) If $\lambda > 1$, then there is a constant $\gamma > 0$ such
that
 $\lim_{N\rightarrow\infty} P \( \sigma_N \geq e^{\gamma N} \) = 1.$
 \end{thm}

\noindent{\it Remarks}: It is not difficult to establish estimates
of the opposite direction, see Theorems \ref{shang} and
\ref{shang1}.  Together we have shown that $\sigma_N$ increases
logarithmically if $\lambda$ is small enough and exponentially if
$\lambda > 1$.

 For any non-empty set $A =\{x_1, x_2 \cdots,
x_k\}$, we assume without loss of generality that  $x_1 < x_2 <
\cdots < x_k$ and define
 $$
 \nu_\psi(A) = \left\{ \begin{array}{ll}
\psi(x_2-x_1) \psi (x_3-x_2) \cdots \psi(x_k-x_{k-1}) & {\rm if} \
k > 1;
 \\
1 & {\rm if} \ k = 1.
 \end{array} \right.
 $$
Let 
 $
{\cal S}_N = \{ 0,1 \}^{\{ 1, \cdots, N \}}$, $K_N = \sum_{A\in
{\cal S}_N \setminus \{ \emptyset \}} \nu_\psi (A)$,  {\rm and}
$\pi (A) =\nu_\psi(A) /K_N.
 $
Then $\pi$ is a probability measure on ${\cal S}_N$.

 \begin{thm}\label{cth}
Suppose that $\lambda = 1$ and the initial distribution  is $\pi$.
Then
 $$
\lim_{N \rightarrow + \infty} P \( \frac N  {C_N}  \leq \sigma_N
\leq C_N N^2 \) = 1.
 $$
 \end{thm}

We now proceed to prove Theorems \ref{shang2} and \ref{cth} by
three different approaches.


\section{Comparison by Coupling }


We will prove the first part of Theorem \ref{shang2}  by
establishing a more general conclusion (Theorem \ref{shang3}). Let
$\{X_t : t \ge 0\}$ be a birth and death process on $\{ 0, 1,
\cdots, N \}$ with
 $$
\begin{array}{llll}
\mbox{ death rate}:  \ & a_i  = & i, \quad &\mbox{ for } i = 1,
\cdots,
N; \\
\mbox{ birth rate}: \ & b_i  = & (i+1)  \alpha ,\quad & \mbox{ for
} i = 0, \cdots, N-1.
 \end{array}
 $$
Let $\tau  = \inf \{ t > 0 : X_t = 0 \}$ be the first time that
$\{X_t : t \ge 0\}$ hits $0$. Let $E^N$ be the conditional
expectation  on $X_0= N$.

 \begin{lem}\label{ex}
 Suppose that $X_0= N$. For
large $N$, 
 $$
  E^N \tau \leq \left\{ \begin{array}{ll}
(2\log N)/(1-\alpha) & {\rm if} \  \alpha < 1;
 \\
2N \log N & {\rm if} \  \alpha = 1;
 \\
\alpha^N \alpha /(\alpha-1)^2 & {\rm if} \ \alpha > 1.
 \end{array} \right.
 $$
 Furthermore
 \begin{equation} \label{zia}
E^N  \tau^2  \leq  2 \( E^N \tau \)^2.
 \end{equation}
 \end{lem}

\noindent{\bf Proof}. Let $P^i$ be the conditional probability
distribution on the initial state $i$, $E^i$ be the expectation
with respect to $P^i$, and $m_i = E^i \tau$ for $i = 0, \cdots,
N$. It is shown in Wang (1980) that
 $$
E^N \tau = \sum_{i=1}^N e_i, \ \ \ E^N  \tau^2 = \sum_{i=1}^N
\varepsilon_i,
 $$
where
 \begin{eqnarray}
e_i
 & = &
\frac{1}{a_i}+ \sum_{k=0}^{N-1-i}\frac{b_ib_{i+1} \cdots
b_{i+k}}{a_ia_{i+1} \cdots a_{i+k}a_{i+k+1}}\\
 & = & \frac 1 i (1  +  \alpha +  \alpha^2 + \cdots  \alpha^{N-i}). \label{E.ei} \\
 \varepsilon_i
 & = &
\frac{2 m_i}{a_i} + \sum_{k=0}^{N-1-i} \frac{2 b_i b_{i+1} \cdots
b_{i+k} m_{i+k+1} } {a_i a_{i+1} \cdots a_{i+k}
a_{i+k+1}}\nonumber
 \end{eqnarray}
 Notice that $ m_i \leq
m_N $ for any $i \leq N$. It follows that $ \varepsilon_i \le 2
m_N e_i$. Therefore,
 $$
E^N  \tau^2  = \sum_{i=1}^N \varepsilon_i \le 2 m_N \sum_{i=1}^N
e_i \leq 2 m_N E^N {\tau} = 2 \( E^N \tau \)^2.
 $$

If $\alpha= 1$,  by (\ref{E.ei}), $e_i = (N-i+1)/i$, and for large
$N$,
$$
E^N \tau = \sum_{i=1}^N e_i \le N \sum_{i=1}^N i^{-1} \le 2 N \log
N. $$

If $\alpha < 1$, 
 $
E^N \tau = \sum_{i=1}^N e_i \le ( 1- \alpha)^{-1} \sum_{i=1}^N
i^{-1} \leq (2 \log N)/( 1- \alpha ); $

If $\alpha > 1$, 
then
 $
E^N \tau = \sum_{i=1}^N e_i \le (  \alpha -1 )^{-1} \sum_{i=1}^N
 \alpha^{N-i+1} \leq   \alpha^{N+1} /(  \alpha -1 )^2$.
 \qed

Consider a nearest particle system $\{ \xi^N_t : t \ge 0 \}$ on
$\{1,2,\cdots,N\}$  starting from
  $\{1,2,\cdots, N\}$ (not necessarily reversible).  Let $\sigma_N$
be the first hitting time of the empty set by $\xi^N_t$, and
$$
M = \max \{\max_n \sum_{l+r = n} \beta(l, r),\quad \sum_l \beta(l,
\infty)\}.$$

 \begin{thm} \label{shang3} Suppose the initial state is $\{1,2,
 \cdots, N\}$. If $M< 1$, then $E \sigma_N \leq (2\log N)/(1-M)$; and
for any sequence $\{ C_N : N \ge 1 \}$ of increasing numbers such
that $\lim_{N \rightarrow \infty} C_N = \infty$,  $ \lim_{N
\rightarrow \infty} P \big( \sigma_N \leq C_N log N\big) = 1$.
 \end{thm}
\medskip
\noindent  {\bf Proof.}  Let $|A|$ be the cardinality of set $A$.
For any configuration $\xi$ such that $|\xi| = i$, there are at
most $i+1$ intervals of consecutive vacant sites, separated by
occupied sites; the rate that a new particle in each interval is
born is no more than $ M$. Hence the rate that $|\xi^N_t|$
increases by 1 is no more than $(i+1) M$. On the other hand, when
$|\xi_t| = i$, the rate that $| \xi_t |$ decreases by 1 is equal
to $i$, the total number of particles.  Compare $|\xi_t|$ with a
birth and death process $X_t$ with parament $\alpha = M$. Since
initially $X_0 = |\xi^N|$, there is a coupling of $\{ X_t : t \ge
0 \}$ and $\{ \xi^N_t : t \ge 0 \}$ such that
 \begin{equation} \label{control}
P^{N, \xi^N} \( X_t \geq |\xi^N_t|, \ \forall \ t \geq 0 \) = 1,
 \end{equation}
where $P^{N,\xi^N}$ is the coupling measure with the initial state
$(N,\xi^N)$. By (\ref{control}), $\sigma_N$ is stochastically
dominated by $\tau$, 
$i.e$., for any $t \ge 0$,
 \begin{equation} \label{E.estimate1}
P (\sigma_N > t) \le P^N (\tau \geq t).
  \end{equation}
By the Chebyshev inequality and  (\ref{zia}), for any $c_N > 0$,
 \begin{equation} \label{E.estimate2}
P (\sigma_N > c_N E^N \tau) \le  P^N \left( \tau \geq c_N E^N \tau
\right)
 \le \frac
{E^N \tau^2 } {\(c_N E^N \tau \)^2} \le \frac 2 {c_N^2}.
  \end{equation}
For any sequence $C_N \rightarrow \infty$ as $N \rightarrow
\infty$, choose $c_N = C_N (1-M)/2$. Then an upper estimate of
$\sigma_N$ may be taken as $c_N E^N \tau$, and the claims in
Theorem \ref{shang3} hold by (\ref{E.estimate2}) and Lemma 2.1.
\qed

By the same argument it is not difficult to establish following
estimates, though a renormalization argument is used in the proof
of the second part of Theorem \ref{shang1}. We will skip the
proof, since they are not needed in proving Theorems \ref{shang2}
and \ref{cth}.

 \begin{thm} \label{shang} Suppose the initial state is $\{1,2,
 \cdots, N\}$. \\
(1) If $M =1$,  then $E \sigma_N \leq N \log N$;  and
 $
\lim_{N \rightarrow \infty} P \big( \sigma_N \leq C_N N log N\big)
= 1; $\\
 (2) If $ M >1$,  then $E \sigma_N \leq M^{N+1}/(M-1)^2$;  and there is a constant $\gamma_1 > 0$ such that
 $$
\lim_{N\rightarrow\infty} P \( \sigma_N \leq e^{\gamma_1 N} \) =
1.
 $$
 \end{thm}

 \begin{thm} \label{shang1} Suppose the initial state is $\{1,2,
 \cdots, N\}$.\\
(1) For any $ \varepsilon > 0$,  $ \lim_{N \rightarrow \infty} P
\big(  \sigma_N > (1-\varepsilon) \log N \big) = 1$;\\
 (2) If $\max_n \min\{ \frac 1 2 \sum_{l = n}^{2n} \beta(l, 3n-l),
\quad \sum_{l=n}^{2n} \beta(l, \infty)\}$ 
is larger than the critical value of the contact process on $\Z$,
then there is a constant $\gamma > 0$ such that
 $$
\lim_{N\rightarrow\infty} P \( \sigma_N \geq e^{\gamma N} \) = 1.
 $$
 \end{thm}

\section{A Lower Estimate of $\sigma_N$}


We first extend the notation  introduced before Theorem \ref{cth}.
For any non-empty set $A =\{x_1, x_2 \cdots, x_k\}$,  $x_1 < x_2 <
\cdots < x_k$,  define
 $$
 \nu_{\psi,\lambda}(A) = \left\{ \begin{array}{ll}
\lambda^{k-1}\psi(x_2-x_1) \psi (x_3-x_2) \cdots \psi(x_k-x_{k-1})
& {\rm if} \ k > 1;
 \\
1 & {\rm if} \ k = 1.
 \end{array} \right.
 $$
Let $ {\cal S}_N = \{ 0,1 \}^{\{ 1, \cdots, N \}}$, $K_N (\lambda)
= \sum_{A\in {\cal S}_N \setminus \{ \emptyset \}}
\nu_{\psi,\lambda} (A)$, {\rm and} $\pi (A) =\nu_{\psi,\lambda}(A)
/K_N(\lambda)$.  Then $\pi$ is a probability measure on ${\cal
S}_N$.

\begin{lem}\label{KN}    $K_N  (\lambda) \ge C N^2 e^{\gamma(\lambda)
N}$ for $\lambda  \geq 1$, where $\gamma (1) = 0$ and $\gamma
(\lambda) > 0$ if $\lambda > 1$.
\end{lem}

\medskip\noindent {\bf Proof.}
 \begin{equation} \label{E.Kn}
K_N (\lambda) = \sum_{\xi \in \mathcal{S}_N \setminus \{ \emptyset
\}} \nu_{\psi,\lambda} (\xi) \ge \sum_{x=0}^{[N/3]} \sum_{y= [2N /
3]}^N \sum_{\xi \in S_N (x, y)}\lambda^{|\xi|-1} \nu_\psi (\xi),
 \end{equation}
where
$$ S_N (x, y) = \left\{ \xi\in {\cal S}_N :  \xi (x) = \xi
(y) = 1, \xi (z)=0, \ \forall \ 1\leq z < x, {\rm \ or \ } y< z
\leq N \right\}.
 $$
In light of (\ref{momn}), by the Renewal Theorem,
  $\nu_\psi (S_N (x,y)) \ge 1/(2\sum_n n\psi(n))$ whenever $y-x$ is large enough.
  If $\lambda =1$, then $K_N  \ge C N^2$ when $N$ is large, and we are done.
If $\lambda > 1$, we can choose  constant  $\delta >0$ such that
$$\nu_\psi \big(\{\xi\in S_N (x,y); |\xi| \geq \delta |y-x|\}\big)\geq
\frac {\nu_\psi (S_N (x,y))} 2.$$
 This together with (\ref{E.Kn}) implies the desired conclusion.
 In particular we may choose $\gamma(\lambda) = (\delta/3)\log\lambda$. \qed

We now use an idea in proving Theorem 7.1.20 of Liggett (1985) to
prove
 \begin{equation} \label{gg2}
\lim_{N\to\infty} P^\pi  \( \sigma_N \ge \frac {K_N}{C_N N} \) =
1.
 \end{equation}
The  first half of Theorem \ref{cth} readily follows from
(\ref{gg2}) and Lemma \ref{KN}. Notice that the hitting time of
the nearest particle system starting from $\{1,2,\cdots, N\}$ is
stochastically larger than that starting from the initial
distribution $\pi$. Therefore the second part of Theorem
\ref{shang2} also follows, with a little change in $\gamma$.

\medskip\noindent {\bf Proof of  (\ref{gg2}).}
The reversible nearest particle system $\{ \xi^N_t : t \ge 0 \}$
is a Markov process taking values in $\mathcal{S}_N$ with jump
rate
$$
q (A, B) = \left\{
 \begin{array}{ll}
1 & {\rm if} \ x\in A,\ B = A\setminus \{x\};
 \\
\lambda \frac{\psi(l_x (A))\psi(r_x (A))}{\psi(l_x (A)+ r_x (A))}
& {\rm if} \ x \notin A,\ B = A \cup \{x\};
 \\
0 & {\rm otherwise.}
 \end{array} \right.
 $$
It is reversible with respect to $\pi$ in the sense that
$\pi(A)q(A, B) = \pi (B) q(B, A)$ for $A, B \in \mathcal{S}_N
\setminus \{ \emptyset \}$.

Let $\{ \widetilde{\xi^N_t} : t \ge 0 \}$ be a Markov process on
$\mathcal{S}_N$, which is a modification of $\{ \xi^N_t : t \ge 0
\}$ so that particles can be born from the empty set. More
specifically,  the transition rates of $\{ \widetilde{\xi^N_t} : t
\ge 0 \}$ is defined as follows.
 $$
\tilde{q} (A, B) = \left\{
 \begin{array}{ll}
q (A, B) & {\rm if} \ A \neq \emptyset;
 \\
q & {\rm if} \ A = \emptyset \ {\rm and} \ |B| = 1;
 \\
0 & {\rm otherwise,}
 \end{array} \right.
 $$
where $q > 0$ is a constant to be determined later. Let $K_N$
stand for $K_N(\lambda)$,
 $$
\nu_\psi \( \{\emptyset\} \) = q^{-1}, \ \ \mbox{and } \ \  \
\tilde\pi = \nu_\psi / \( K_N + q^{-1} \).$$
 Then $\{ \widetilde{\xi^N_t} : t \ge 0 \}$ is reversible with
respect to $\tilde \pi$ in the sense that $\tilde \pi(A) q(A, B) =
\tilde \pi (B) q(B, A)$ for any $A, B\in \mathcal{S}_N$.

 Let $\tilde{P}$ be the distribution of $\{ \widetilde{\xi^N_t} : t \ge
0 \}$ with initial distribution $\tilde \pi$, and $\tilde{E}$ be
the expectation with respect to $\tilde{P}$. Notice that $\{
\widetilde{\xi^N_t} : t \ge 0 \}$ is stationary under $\tilde{P}$.
For any $t > 0$,
 $$
2 t \tilde{\pi} (\{\emptyset\})
 =
\tilde{E} \int_0^{2t} 1_{ \{ \widetilde{\xi^N_s} = \emptyset \} }
ds.
 $$
 Introduce the stopping time
 $
\tau = \inf \{ t \geq 0 : \widetilde{\xi^N_t} = \emptyset \}.
 $
By the Strong Markovian Property, the right side above equals
 \begin{eqnarray*}
 & &
\tilde{E}  \tilde{E} \( \left. \int_0^{2t}
1_{\{\widetilde{\xi^N_s} = \emptyset\}} d s \right|
\mathcal{F}_\tau \)
 \ge
\tilde{E} \tilde{E} \(  \left. 1_{\{\tau<t\}} \int_0^{2t}
1_{\{\widetilde{\xi^N_s} = \emptyset\}} d s \right|
\mathcal{F}_\tau \)
 \\ & \ge &
\tilde{E} \tilde{E} \( 1_{\{ \tau < t \}}  \left. \int_\tau^{\tau
+ t} 1_{\{\widetilde{\xi^N_s} = \emptyset\}} d s \right|
\mathcal{F}_\tau \)
 =
\tilde{P} (\tau<t) \tilde{E} \( \left. \int_0^t
1_{\{\widetilde{\xi^N_s} = \emptyset\}} d s \right|
\widetilde{\xi^N_0} = \emptyset \)
 \end{eqnarray*}
Denote by $\sigma$ the first time $\{\widetilde{\xi^N_t} : t \ge 0
\}$ jumps. Then
 \begin{eqnarray*}
\tilde{E} \( \left. \int_0^t 1_{\{\widetilde{\xi^N_s} =
\emptyset\}} d s \right| \widetilde{\xi^N_0} = \emptyset \)
  & \ge &
\tilde{E} \( \sigma 1_{\{ \sigma \le t \}} | \widetilde{\xi^N_0} =
\emptyset \)
 =
\int_{0}^{t} s \tilde{q}_\emptyset e^{ - \tilde{q}_\emptyset s } d
s,
 \end{eqnarray*}
where $\tilde{q}_\emptyset = \sum_{\xi} \tilde{q} (\emptyset, \xi)
= N q$. Hence
 \begin{equation} \label{2t}
\tilde{P} (\tau<t) \le\frac { 2t \tilde{\pi} (\{\emptyset\})}{
\int_{0}^{t} s \tilde{q}_\emptyset e^{ - \tilde{q}_\emptyset s } d
s } = \frac{ 2 t q^{-1}}{ K_N + q^{-1}} \cdot \frac{Nq}{
1-e^{-Nqt} - Nqt e^{-Nqt}}.
 \end{equation}
On the other hand,
 \begin{eqnarray*}
\tilde{P} (\tau<t)
 & \geq &
\tilde{P} \( \tau<t,\  \widetilde{\xi^N_0} \neq \emptyset \)
 =
\tilde{P} \( \widetilde{\xi^N_0} \neq\emptyset \) \tilde{P} \(
\tau<t|\widetilde{\xi^N_0}\neq\emptyset \)
  \\ & = &
\frac {K_N}{ K_N + q^{-1} } P (\sigma_N<t) .
 \end{eqnarray*}
This together with $(\ref{2t})$ yields that
 $$
P \( \sigma_N < t \) \le \frac {2 t N} {K_N \( 1-e^{-Nqt} - Nqt
e^{-Nqt} \)} \ .
 $$
Let $q \rightarrow \infty$, then
 $$
P \( \sigma_N < t \) \le  \frac {2 t N } {K_N} \ .
 $$
This implies (\ref{gg2}), by choosing $t = K_N/(C_N N)$. \qed


\section{The Critical Case}

In this section we will prove the second half of Theorem
\ref{cth}, $i.e.$,  when $\lambda =1$,
 \begin{equation} \label{gg1}
\lim_{N\to\infty} P^\pi \( \sigma_N \le C_N N^2 \) = 1.
 \end{equation}
Let $\{\eta_t : t \ge 0 \}$ be an infinite reversible nearest
particle system on $\Z$ with finite many particles to the right of
the origin (The third case on page 2); and $r_t$ the rightmost
particle in $\{ \eta_t : t \ge 0 \}$, i.e.
 $r_t : = \sup \{ x : \eta_t (x) = 1 \}$.
The properties of $r_t$ of the critical nearest particle system
are studied in Schinazi (1992). For a recent survey, see Mountford
(2003).

 \begin{lem}\label{edge} {\rm (Schinazi (1992), Theorem 1)}
 Let $\{\eta_t : t \ge 0 \}$
be the critical reversible nearest particle system on $\Z$.
Suppose the initial configurations have a particle at the origin
and no particle to  the right of the origin, and follows the
renewal measure $Ren(\psi)$ with density $\psi(\cdot)$. Then, as
$a\rightarrow\infty$, $r_{a^2t}/a$ converges in distribution to a
Brownian motion with diffusion constant $D > 0$ in the Skorohod
space.
 \end{lem}

\noindent {\bf Proof of (\ref{gg1}).} 
Partition the configuration space ${\cal S}_N$ according to the
position of the rightmost particle. Namely, let
 $$
A_x = \left\{ \xi \in {\cal S}_N : \xi (x) = 1, \mbox{and}\ \xi(y)
= 0 \ \mbox{ for any }\ y > x \right\}
 $$
be the set of configurations whose rightmost particle is at $x$.
Denote by $P$ the distribution of $\{ \xi^N_t : t \ge 0 \}$ with
initial distribution $\pi$, and by $P_{N,x}$ the conditional
distribution of the nearest particle system on $\{1,2, \cdots,
N\}$ whose initial configurations are in $A_x$.  Then
 \begin{equation}\label{E.PartitionP}
P = \sum_{x=0}^N P(A_x) P_{N,x}.
 \end{equation}

Denote by $\P$ the distribution of the nearest particle system on
$\Z$ with the initial distribution in Lemma \ref{edge}, and $\P_x$
the translation of $\P$ by $x$. Thanks to the attractive property,
there is a coupling of $\P_x$ and $P_{N,x}$ such that for all $t
> 0$ and all $i\in \Z$,
 \begin{equation}\label{xiao}
\xi^N_t (i) \le \eta_t (i).
 \end{equation}
Then under this coupling, $\xi^N_t \equiv \emptyset$ once $r_t <
1$, hence $ \sigma_N \le \inf \{t : r_t < 1 \}$.

Suppose that $\lim_{N \rightarrow \infty} C_N = \infty$. For any
$C > 0$ and large $N$,
  \begin{eqnarray*}
P_{N,x} \( \sigma_N \le C_N N^2 \) & \ge &P_{N,x} \( \sigma_N \le
C (x-1)^2 \) \\
 & \ge &
\P_x \( \exists \ t \le C (x-1)^2 {\rm \ s.t.} \ r_t < 1 \)\\
 & = &
\P \( \exists \ t \le C (x-1)^2 {\rm \ s.t.} \ r_t < - (x-1) \)
 \\ & = &
\P \( \exists \ t \le C {\rm \ s.t.} \ r_{ (x-1)^2 t} / (x-1) < -
1 \).
 \end{eqnarray*}
Here the first equality holds because $\P_x$ is the translation of
$\P$ by $x$.  This together with Lemma \ref{edge} implies that
 $$
\liminf_{N, x \rightarrow + \infty} P_{N,x} \( \sigma_N \le C_N
N^2 \) \ge \P \( \exists \ t \le C {\rm \ s.t.} \ B_t < - 1 \), \
\ \ \forall \ C > 0,
 $$
where $\{ B_t : t \ge 0 \}$ is a Brownian motion with diffusion
constant $D$. Let $C \rightarrow + \infty$, the right side of the
above equation converges to 1. Hence
 $$
\lim_{N,x \rightarrow + \infty} P_{N,x} \( \sigma_N \le C_N N^2 \)
= 1.
 $$
Consequently, for any $\varepsilon > 0$, there exists $N_0
> 0$ such that for any $N \ge x \ge N_0$
 $$
P_{N,x}  \( \sigma_N \le C_N N^2 \) > 1 - \varepsilon.
 $$
This together with (\ref{E.PartitionP}) implies that
 \begin{equation} \label{E.P}
P  \( \sigma_N \le C_N N^2 \) = \sum_{x=1}^N P(A_x) P_{N,x}  \(
\sigma_N \le C_N N^2 \) \ge (1 - \varepsilon) \sum_{x = N_0}^N
P(A_x).
 \end{equation}
On the other hand,
 $$
\sum_{x = 1}^{ N_0 - 1} \nu_\psi (A_x) \le \sum_{x = 1}^{ N_0 - 1}
\sum_{y = 1 }^x \nu_\psi (S_N (y,x)) \le N_0^2.
 $$
Therefore, as $N \rightarrow \infty$,
 $$
\sum_{x = N_0}^N P(A_x) \ge 1 - N_0^2 / (C  N^2) \rightarrow 1.
 $$
This together with (\ref{E.P}) implies that
 $
\liminf_{N \rightarrow \infty} P  \( \sigma_N \le C_N N^2 \) \ge 1
- \varepsilon.
 $
Let $\varepsilon \rightarrow 0$ and the result follows.\qed

\vspace{0.5cm}
\small

\baselineskip=0.7\baselineskip

\noindent {\bf  References}

\noindent Durrett, R. and Liu, X. F. (1988). The contact process
on a finite set. {\em Ann. Probab.} {\bf 16} 1158--1173.

\noindent Durrett, R. and Schonmann, R. H. (1988). The contact
process on a finite set II. {\em Ann. Probab.} {\bf 16}
1570--1583.

\noindent Durrett, R., Schonmann, R. H. and Tanaka, N. I. (1989).
The contact process on a finite set III: The critical case. {\em
Ann. Probab.} {\bf 17} 1303--1321.

\noindent Liggett, T. M. (1985). {\em Interacting particle
systems.} New York, Springer-Verlag.

\noindent Mountford, T.S. (1992). A critical value for the uniform
nearest particle system,  {\em Ann. Probab.} {\bf 20} 2031--2042.

\noindent Mountford, T.S. (2003). Critical reversible attractive
nearest particle systems, In {\em Topics in Spatial Stochastic
Processes, Lecture Notes in Mathematics \bf 1802}, Springer,
Berlin.

\noindent Schinazi, R. (1992). Brownian fluctuations of the edge
for critical reversible nearest particle systems. {\em Ann.
Probab.} {\bf 20} 194--205.

\noindent Wang Z. K. (1980). {\em Birth and Death Processes and
Markov Chains} (in Chinese). Beijing, Science Publishing House.

\vspace{0.5cm} \noindent LMAM, School of Mathematical Sciences,
Peking University, Beijing 100871, China
\end{document}